\documentclass[11pt]{article}

\usepackage[margin=1in]{geometry}
\usepackage{amsmath,amssymb,amsthm,mathtools}
\usepackage{enumitem}
\usepackage[colorlinks=true,linkcolor=blue,citecolor=blue,urlcolor=blue]{hyperref}

\numberwithin{equation}{section}

\newtheorem{theorem}{Theorem}[section]
\newtheorem{lemma}[theorem]{Lemma}
\newtheorem{proposition}[theorem]{Proposition}
\newtheorem{corollary}[theorem]{Corollary}
\newtheorem{definition}[theorem]{Definition}
\newtheorem{assumption}[theorem]{Assumption}
\newtheorem{remark}[theorem]{Remark}

\newcommand{\Om}{\Omega}
\newcommand{\R}{\mathbb R}
\newcommand{\C}{\mathbb C}
\newcommand{\supp}{\operatorname{supp}}
\newcommand{\dist}{\operatorname{dist}}
\newcommand{\dH}{d_{\mathrm H}}
\newcommand{\calN}{\mathcal N}
\newcommand{\pnu}{\partial_\nu}
\newcommand{\eps}{\varepsilon}

\DeclareMathOperator{\diam}{diam}

\title{Localization of Sparse Passive Sources with Unknown Wave Speed\\
from Low-Frequency Boundary Data}
\author{Amir Moradifam\footnote{Department of Mathematics, University of California, Riverside, California, USA. E-mail: amirm@ucr.edu}}

\date{}

\begin{document}

\maketitle

\begin{abstract}
Simultaneously determining an initial source and an unknown wave speed from a single boundary measurement remains open in general, while the unrestricted linearized problem is unstable in every Sobolev scale. We show that spatial sparsity overcomes this obstruction at the level of source localization. For initial sources given by finite sums of Dirac masses, the source locations and the weighted amplitudes $\alpha_jc(x_j)^{-2}$ can be determined uniquely and stably from low-frequency boundary data, without any knowledge of the wave speed. These results require no smoothness or smallness assumptions on the wave speed, only uniform positive upper and lower bounds. We further extend the localization result to sparse sources whose individual components have small spatial support and may have different shapes and sizes. This moves the theory beyond the idealized point-source model and allows for more realistic localized sources.
\end{abstract}

\section{Introduction}\label{sec:introduction}

Inverse source problems for wave equations arise when waves generated by
unknown internal events are observed on the boundary of a region and one seeks
to infer properties of the source. In many passive imaging problems, the
primary objective is not to reconstruct an arbitrary source profile but to
locate a small number of spatially concentrated events. This is a useful
idealization in microseismic monitoring, acoustic-emission testing, structural
health monitoring, passive nondestructive evaluation, and imaging of small
photoacoustic absorbers
\cite{AkiRichards2002,GesretEtAl2015,GrosseOhtsu2008,HassanEtAl2021}.
In such settings, the propagation medium may be imperfectly calibrated, while
the source locations remain the quantities of principal interest.

We consider the scalar wave equation
\begin{equation}\label{eq:intro-wave}
\begin{cases}
c(x)^{-2}\partial_t^2u-\Delta u=0,
    & (t,x)\in(0,\infty)\times\R^3,\\
u(0,x)=f(x),
    & x\in\R^3,\\
\partial_tu(0,x)=0,
    & x\in\R^3,
\end{cases}
\end{equation}
where \(c(x)>0\) is the wave speed and \(f\) is the initial source. Let
\(\Om\Subset\R^3\) contain both the source and the inhomogeneity of the
medium. The measured quantity is the boundary trace
\[
\Lambda_{f,c}
=
u|_{(0,\infty)\times\partial\Om},
\]
understood distributionally when \(f\) is a measure. The negative-Sobolev
formulation of the Cauchy problem and the distributional boundary restriction
are described in Section~\ref{sec:low-frequency}; the low-frequency
formulation may also be read directly in the elliptic Laplace domain. Our
inverse problem is to determine the support of a sparse source without
assuming that \(c\) is known.

At first sight, localization without a calibrated speed appears problematic.
Changes in the propagation speed alter travel times and can therefore affect
locations inferred from time-domain propagation. More fundamentally,
simultaneous recovery of an unrestricted source and an unknown speed from a
single boundary measurement is considerably more delicate than either
inverse problem separately. The main observation of this paper is that finite
spatial sparsity removes this ambiguity at zero frequency. Although the speed
remains unknown, it enters the leading low-frequency coefficient only through
nonzero weights attached to the source points. The source locations can
therefore be recovered without first determining the speed.

When the wave speed is known, recovery of an initial source in photoacoustic
and thermoacoustic tomography is well understood under appropriate geometric
hypotheses; see, among many references,
\cite{AgranovskyKuchmentKunyansky2009,AmmariBossyJugnonKang2010,
Hristova2009,HristovaKuchmentNguyen2008,KuchmentKunyansky2008,
StefanovUhlmann2009,StefanovUhlmann2011,StefanovUhlmann2013SourceSpeed,
Wang2009,XuWang2006}. The effect of uncertainty in the sound speed has also
been studied in reconstruction and stability analysis
\cite{HristovaKuchmentNguyen2008,OksanenUhlmann2014}. By contrast,
positive results for the simultaneous recovery of a source and an unknown
speed require additional structure
\cite{FinchHickmann2013,KianUhlmann2025,KnoxMoradifam2020,LiuUhlmann2015},
and the unrestricted linearized problem is unstable in every Sobolev scale
\cite{StefanovUhlmann2013Instability}.
More recently, Kian and Liu~\cite{KianLiu2026} established uniqueness and
H\"older stability for simultaneous recovery under piecewise-constant wave-speed
assumptions and a quantitative closeness condition relative to a constant background.

The additional structure used here is sparsity. Suppose that
\begin{equation}\label{eq:intro-atomic}
f
=
\sum_{j=1}^{N}\alpha_j\delta_{x_j},
\qquad
\alpha_j\neq0.
\end{equation}
Recovery of finitely many point sources from static potential data, or from
Cauchy data for elliptic equations, has a substantial literature. Qualitative
and quantitative results appear, for example, in
\cite{ElBadiaHaDuong2000,Vessella1992,ElBadiaElHajj2012,
KandaswamyBluVanDeVille2009,ElBadiaNara2011,BaoLiuTriki2021}.
Small-absorber and wave-based variants have also been studied
\cite{AmmariBossyJugnonKang2010,AlJebawyElBadia2022,
AlJebawyElBadiaTriki2022}. Our contribution is to connect the passive wave
measurement with unknown speed to an exact static problem for a weighted
source and to develop the resulting uniqueness and stability consequences
for atomic and mollified sources.

Write
\[
q(x)=c(x)^{-2}.
\]
Motivated by the low-frequency method used in
\cite{KnoxMoradifam2020}, we apply a one-sided Laplace transform to
\eqref{eq:intro-wave}. After division by the transform parameter and passage
to zero frequency, the transformed field converges to the decaying solution
of
\begin{equation}\label{eq:intro-poisson}
-\Delta W=qf
\qquad\text{in }\R^3.
\end{equation}
Thus the leading zero-frequency boundary coefficient is the trace of the
Newtonian potential generated by \(qf\). For two source-speed pairs, the
discrepancy between these coefficients controls all harmonic moments of
\[
\mu=q_2f_2-q_1f_1.
\]
If \(f\) is given by \eqref{eq:intro-atomic}, then
\begin{equation}\label{eq:intro-weighted-atomic}
qf
=
\sum_{j=1}^{N}\alpha_jq(x_j)\delta_{x_j}.
\end{equation}
Since \(q>0\), multiplication by \(q\) changes the amplitudes but not the
support. The unknown-speed localization problem is therefore reduced to the
determination of a finite weighted atomic measure from its harmonic moments.

We first show that harmonic moments uniquely determine a finite atomic
measure. It follows that equality of the time-domain boundary traces for two
source-speed pairs implies
\[
q_1f_1=q_2f_2
\qquad\text{and}\qquad
\supp f_1=\supp f_2.
\]
The data therefore determine not only the source locations but also the
weighted amplitudes
\[
\alpha_jq(x_j)
=
\alpha_jc(x_j)^{-2}.
\]
The same reduction also yields quantitative stability. We consider classes
with a fixed upper bound on the number of atoms, a prescribed lower bound on
their separation, and a positive lower bound on the nonzero weighted
amplitudes \( |\alpha_j c(x_j)^{-2}| \). If
\[
\eps
=
\|\mathcal L_{f_2,c_2}-\mathcal L_{f_1,c_1}\|_
 {H^{1/2}(\partial\Om)}
\]
denotes the discrepancy between the corresponding zero-frequency boundary
coefficients, we first obtain a global conditional H\"older estimate. For
sufficiently small \(\eps\), this estimate improves to
\begin{equation}\label{eq:intro-lipschitz}
\dH(\supp f_1,\supp f_2)
\le C\eps.
\end{equation}
In this regime, the separation assumption also gives a unique matching of the
source points. Thus, if
\[
f_1=\sum_{j=1}^{N}\alpha_j\delta_{x_j},
\qquad
f_2=\sum_{j=1}^{N}\beta_j\delta_{y_j},
\]
then the points may be relabeled so that
\[
|x_j-y_j|\le C\eps.
\]
Under an additional uniform upper bound on the weighted amplitudes, we prove
the joint estimate
\begin{equation}\label{eq:intro-joint-stability}
\max_{1\le j\le N}
\left(
|x_j-y_j|
+
|\alpha_jq_1(x_j)-\beta_jq_2(y_j)|
\right)
\le C\eps.
\end{equation}
Hence the complete weighted atomic measure \(qf\), which is the natural
identifiable object in this problem, is locally Lipschitz stable.
Notably, apart from the standing assumption that \(c=1\) outside a compact
subset of \(\Om\), these Lipschitz stability results require no smoothness or
smallness assumptions on the wave speed, only uniform positive upper and lower bounds.

The proof uses complex harmonic polynomials to separate and interpolate
finite point configurations. If \(\zeta\in\C^n\) satisfies
\(\zeta\cdot\zeta=0\), then every polynomial in \(\zeta\cdot x\) is
harmonic. A suitable choice of \(\zeta\) reduces the exact moment identities
to a Vandermonde system and hence gives uniqueness. For stability, harmonic
Lagrange polynomials isolate one atom while vanishing at all other atoms in
the two competing configurations. A global interpolation estimate initially
gives only a H\"older bound with exponent \(1/N_0\). This places the two
configurations within a fixed fraction of the prescribed separation scale.
In that local regime, only one cross-configuration factor can be small, and a
second interpolation estimate gives the Lipschitz bounds
\eqref{eq:intro-lipschitz} and \eqref{eq:intro-joint-stability}.

We also consider sources obtained by replacing the point masses by small,
nonnegative, normalized profiles that may have different shapes and radii.
The atomic assumption is essential to exact determination from harmonic
moments: every function of the form \(\Delta\psi\), with
\(\psi\in C_c^\infty(\Om)\), annihilates all harmonic test functions.
Nevertheless, within a uniformly separated blob class, the translation
parameters remain stable. For nonradial profiles, these parameters need not
coincide with a centroid or a point of maximum amplitude. If \(X_1\) and
\(X_2\) denote the two sets of translation parameters and \(\eta_*\) is the
largest blob radius, we prove
\begin{equation}\label{eq:intro-mollified}
\dH(X_1,X_2)
\le C(\eps+\eta_*).
\end{equation}
The same estimate holds for the supports of the two mollified sources. The
term \(\eta_*\) represents the natural spatial resolution of this source
model.

Two limitations should be kept in view. First, the low-frequency coefficient
determines the weighted amplitudes \(\alpha_jc(x_j)^{-2}\), but does not separately determine
the physical amplitudes \(\alpha_j\) without the speed at the source points.
The locations remain identifiable because multiplication by the positive
factor \(c^{-2}\) does not alter the support. Second, the discrepancy in the
stability estimates is measured in the derived low-frequency coefficient
\(\mathcal L_{f,c}\), which is defined through a Laplace transform over all
positive times. We do not claim a stability estimate in terms of data on a
fixed finite observation interval. Such a result would require additional
propagation, decay, or noise assumptions. The seismic and engineering
examples above are likewise motivations for the scalar model, not claims
about elastic, anisotropic, or attenuating systems.

The paper is organized as follows. Section~\ref{sec:low-frequency} derives the
zero-frequency elliptic problem and converts its boundary trace into harmonic
moment information. Section~\ref{sec:harmonic} develops the qualitative and
quantitative interpolation tools. Section~\ref{sec:atomic} proves uniqueness
and conditional stability for atomic sources, including the joint stability
of their locations and weighted amplitudes. Section~\ref{sec:mollified}
treats mollified point sources and the stability of their translation
parameters.

\section{The low-frequency coefficient and harmonic moments}
\label{sec:low-frequency}

We now justify the reduction from the time-domain wave equation to a static
potential problem. A positive Laplace parameter is convenient here because it
avoids a limiting-absorption argument and selects the decaying resolvent
directly.

Throughout this section, \(\Om\Subset\R^3\) is a bounded \(C^\infty\) domain
with connected exterior and \(K\Subset\Om\). The wave speed \(c\) is a fixed
bounded Borel representative, and \(q=c^{-2}\) satisfies
\begin{equation}\label{eq:q-admissible}
q\in L^\infty(\R^3),\qquad
0<q_-\le q(x)\le q_+<\infty\quad\text{for every }x\in\R^3,\qquad
\operatorname{ess\,supp}(q-1)\Subset\Om.
\end{equation}
The source \(f\) is a finite signed or complex Radon measure supported in
\(K\). Because such a measure may charge sets of Lebesgue measure zero, the
choice of a Borel representative of \(q\) on \(\supp f\) is part of the
source-speed pair. The product \(qf\) is the finite measure defined by
\[
\langle qf,\phi\rangle
=
\int_{\R^3}q(x)\phi(x)\,df(x),
\qquad \phi\in C_c^\infty(\R^3).
\]
In particular, if
\[
f=\sum_{j=1}^{N}\alpha_j\delta_{x_j},
\]
then
\[
qf=\sum_{j=1}^{N}\alpha_jq(x_j)\delta_{x_j}.
\]

We use the standard self-adjoint realization of the wave equation with the
rough coefficient \(q\). Let
\[
\mathsf H=L^2(\R^3,q(x)\,dx),
\]
which coincides as a set with \(L^2(\R^3)\), with equivalent norms.
Consider the operator
\[
Av=-q^{-1}\Delta v,
\qquad
\mathcal D(A)=H^2(\R^3).
\]
The bounds in \eqref{eq:q-admissible} imply that \(A\) is nonnegative
and self-adjoint on \(\mathsf H\), as follows from the standard
quadratic-form realization of the weighted Laplacian. Moreover,
the graph norm of \(A\) is equivalent to the \(H^2(\R^3)\) norm.

For \(r\ge0\), let
\[
\mathsf H_r=\mathcal D((I+A)^{r/2}),
\]
with its graph norm, and let \(\mathsf H_{-r}\) be the corresponding dual
space with pivot \(\mathsf H\). Since \(\mathsf H_2=H^2(\R^3)\) with
equivalent norms and \(H^2(\R^3)\hookrightarrow C_b^0(\R^3)\), the measure
\(f\) defines an element of \(\mathsf H_{-2}\) through
\[
\langle f,v\rangle_{\mathsf H_{-2},\mathsf H_2}
:=
\int_{\R^3}q(x)v(x)\,df(x),
\qquad v\in\mathsf H_2,
\]
and
\[
\|f\|_{\mathsf H_{-2}}
\le C_q\|f\|_{\mathrm{TV}}.
\]
This identification agrees with the usual one for function-valued initial
data in the weighted space \(\mathsf H\).

For the coefficient class \eqref{eq:q-admissible}, we understand the Cauchy
problem
\begin{equation}\label{eq:wave-distributional}
\begin{cases}
q(x)\partial_t^2u-\Delta u=0,
    &(t,x)\in(0,\infty)\times\R^3,\\
u(0,\cdot)=f,\\
\partial_tu(0,\cdot)=0,
\end{cases}
\end{equation}
in this self-adjoint sense. The spectral theorem gives the unique solution
\[
u(t)=\cos(tA^{1/2})f,
\]
with
\[
u\in C(\R;\mathsf H_{-2})
\cap C^1(\R;\mathsf H_{-3})
\]
and the global bound
\begin{equation}\label{eq:global-wave-bound}
\sup_{t\in\R}\|u(t)\|_{\mathsf H_{-2}}
\le
\|f\|_{\mathsf H_{-2}}.
\end{equation}

We next define the boundary measurement. By
\eqref{eq:q-admissible}, there is a neighborhood \(\mathcal U\) of
\(\partial\Om\) on which \(q=1\) almost everywhere. The restriction of the
abstract solution to \((0,\infty)\times\mathcal U\) is a distribution: for
\(\psi\in C_c^\infty(\mathcal U)\), its action is defined by the duality
between \(\mathsf H_{-2}\) and \(\mathsf H_2\). Moreover, it satisfies
\[
(\partial_t^2-\Delta)u=0
\qquad\text{in }\mathcal D'((0,\infty)\times\mathcal U).
\]
Indeed, if a spatial test function is supported in \(\mathcal U\), then
\(A\psi=-\Delta\psi\); testing the abstract wave equation against such
functions gives the displayed local equation.
Microlocal elliptic regularity for this constant-coefficient operator
\cite[Theorem~8.3.1]{HormanderI} gives
\[
\operatorname{WF}(u)
\cap T^*((0,\infty)\times\mathcal U)
\subset
\{(t,x;\tau,\xi)\ne0:|\xi|^2=\tau^2\}.
\]
Let
\[
\Sigma=(0,\infty)\times\partial\Om.
\]
Every nonzero conormal covector to \(\Sigma\) has the form
\[
(t,x;0,\lambda\nu(x)),
\qquad \lambda\ne0,
\]
and is noncharacteristic because
\[
|\lambda\nu(x)|^2-0^2=\lambda^2\ne0.
\]
Consequently,
\[
\operatorname{WF}(u)\cap N^*\Sigma=\varnothing.
\]
H\"ormander's distributional restriction theorem
\cite[Corollary~8.2.7]{HormanderI} therefore shows that
\[
\Lambda_{f,c}
:=
u|_{(0,\infty)\times\partial\Om}
\in
\mathcal D'((0,\infty)\times\partial\Om)
\]
is well defined. This is a microlocal distributional restriction; no
ordinary spatial Sobolev trace of \(u(t,\cdot)\) is asserted.

The bound \eqref{eq:global-wave-bound} implies that, for every \(s>0\),
\[
\int_0^\infty e^{-st}\|u(t)\|_{\mathsf H_{-2}}\,dt
\le
\frac{1}{s}\|f\|_{\mathsf H_{-2}}.
\]
Thus the one-sided Laplace transform
\begin{equation}\label{eq:laplace-transform}
U_s
=
\int_0^\infty e^{-st}u(t,\cdot)\,dt
\end{equation}
converges as a Bochner integral in \(\mathsf H_{-2}\). We set
\begin{equation}\label{eq:scaled-laplace}
V_s=\frac{U_s}{s}.
\end{equation}

\begin{proposition}      
\label{prop:laplace-domain}
For every \(s>0\), the scaled Laplace transform \(V_s\) belongs to
\(L^2(\R^3)\) and is the unique \(L^2(\R^3)\) distributional solution of
\begin{equation}\label{eq:modified-helmholtz}
(-\Delta+s^2q)V_s=qf
\qquad\text{in }\R^3.
\end{equation}
\end{proposition}

\begin{proof}
The scalar identity
\[
\int_0^\infty e^{-st}\cos(t\sqrt\lambda)\,dt
=
\frac{s}{s^2+\lambda},
\qquad \lambda\ge0,
\]
and the spectral theorem give
\[
U_s=s(s^2+A)^{-1}f,
\qquad
V_s=(s^2+A)^{-1}f.
\]
The resolvent maps \(\mathsf H_{-2}\) continuously into
\(\mathsf H=L^2(\R^3,q\,dx)\), which is the same set as
\(L^2(\R^3)\) with an equivalent norm. Pairing
\((s^2+A)V_s=f\) with \(\phi\in C_c^\infty(\R^3)\) gives
\[
\int_{\R^3}V_s(-\Delta\phi+s^2q\phi)\,dx
=
\int_{\R^3}q\phi\,df.
\]
This is precisely \eqref{eq:modified-helmholtz} in the distributional
sense.

If \(w\in L^2(\R^3)\) satisfies
\[
(-\Delta+s^2q)w=0,
\]
then \(-\Delta w=-s^2qw\in L^2(\R^3)\), so \(w\in H^2(\R^3)\).
Testing the equation against \(w\) yields
\[
\int_{\R^3}|\nabla w|^2\,dx
+
s^2\int_{\R^3}q(x)|w|^2\,dx
=0.
\]
Since \(q\ge q_->0\), it follows that \(w=0\).
\end{proof}

For \(s\ge0\), set
\[
G_s(x)=\frac{e^{-s|x|}}{4\pi|x|},
\qquad x\ne0.
\]
Then
\[
(-\Delta+s^2)G_s=\delta_0
\]
in \(\mathcal D'(\R^3)\). For a finite Radon measure \(\mu\), write
\[
\calN\mu:=G_0*\mu
\]
for its Newtonian potential. In particular, define
\[
W_{f,c}:=\calN(qf)=G_0*(qf).
\]
Then \(W_{f,c}\in L^2_{\mathrm{loc}}(\R^3)\) and
\[
-\Delta W_{f,c}=qf
\qquad\text{in }\mathcal D'(\R^3).
\]
At every point \(x\in\R^3\) for which
\[
\int_K\frac{1}{|x-y|}\,d|f|(y)<\infty,
\]
the potential admits the pointwise representation
\[
W_{f,c}(x)
=
\frac{1}{4\pi}
\int_K\frac{q(y)}{|x-y|}\,df(y).
\]
In particular, this representation holds for every
\(x\notin\supp f\). On all of \(\R^3\), \(W_{f,c}\) is understood as a
distribution. It is harmonic in \(\R^3\setminus\supp f\) and tends to zero
as \(|x|\to\infty\).

\begin{proposition}[Zero-frequency limit]
\label{prop:zero-frequency}
Let \(V_s\) be as in Proposition~\ref{prop:laplace-domain}. Then, as
\(s\downarrow0\),
\begin{equation}\label{eq:zero-frequency-convergence}
V_s-W_{f,c}\longrightarrow0
\qquad\text{in }H^2_{\mathrm{loc}}(\R^3).
\end{equation}
In particular,
\begin{equation}\label{eq:boundary-zero-frequency}
V_s|_{\partial\Om}
\longrightarrow
W_{f,c}|_{\partial\Om}
\qquad\text{in }C^\infty(\partial\Om).
\end{equation}
\end{proposition}

\begin{proof}
Write
\[
q=1+Q,
\qquad
Q\in L^\infty(\R^3),
\qquad
\operatorname{ess\,supp}Q\Subset\Om.
\]
Equation \eqref{eq:modified-helmholtz} becomes
\[
(-\Delta+s^2)V_s=qf-s^2QV_s.
\]
Since \(V_s\) is the unique \(L^2(\R^3)\) solution, convolution with
\(G_s\) gives
\begin{equation}\label{eq:lippmann-schwinger}
V_s
=
G_s*(qf)-s^2G_s*(QV_s).
\end{equation}

Fix a ball \(B\) containing \(\overline\Om\); it may be enlarged to contain
any prescribed compact subset of \(\R^3\). Since
\[
0<G_s(x)\le G_0(x),
\qquad 0<s\le1,
\]
and \(|x|^{-1}\in L^2_{\mathrm{loc}}(\R^3)\), Minkowski's integral
inequality gives
\[
\begin{aligned}
\|G_s*(qf)\|_{L^2(B)}
&\le
\int_K
\|G_s(\,\cdot-y)\|_{L^2(B)}\,d|qf|(y)\\
&\le
C_B\|qf\|_{\mathrm{TV}},
\end{aligned}
\]
uniformly for \(0<s\le1\). Schur's test also gives
\[
\|G_s*(Qv)\|_{L^2(B)}
\le
C_{B,q}\|v\|_{L^2(B)},
\qquad 0<s\le1,
\]
because \(Q\in L^\infty(\R^3)\) is supported in \(\Om\Subset B\).
Applying these estimates to \eqref{eq:lippmann-schwinger} yields
\[
\|V_s\|_{L^2(B)}
\le
C_B\|qf\|_{\mathrm{TV}}
+
C_{B,q}s^2\|V_s\|_{L^2(B)}.
\]
After decreasing \(s_0>0\), if necessary, the last term can be absorbed
into the left-hand side, and hence
\begin{equation}\label{eq:uniform-local-L2}
\sup_{0<s\le s_0}\|V_s\|_{L^2(B)}
\le
C_B'\|qf\|_{\mathrm{TV}}.
\end{equation}

The inequality \(1-e^{-r}\le r\) gives
\[
|G_s(x)-G_0(x)|
=
\frac{1-e^{-s|x|}}{4\pi|x|}
\le
\frac{s}{4\pi}.
\]
Therefore,
\[
\|(G_s-G_0)*(qf)\|_{L^2(B)}
\le
\frac{s|B|^{1/2}}{4\pi}\|qf\|_{\mathrm{TV}}
\longrightarrow0.
\]
Using \eqref{eq:lippmann-schwinger}, \eqref{eq:uniform-local-L2}, and
Schur's estimate, we obtain
\[
\begin{aligned}
\|V_s-W_{f,c}\|_{L^2(B)}
&\le
\|(G_s-G_0)*(qf)\|_{L^2(B)}
+s^2\|G_s*(QV_s)\|_{L^2(B)}\\
&\longrightarrow0.
\end{aligned}
\]

Set \(z_s=V_s-W_{f,c}\). Subtracting the equations for \(V_s\) and
\(W_{f,c}\) gives
\[
-\Delta z_s=-s^2qV_s
\qquad\text{in }\mathcal D'(\R^3).
\]
If \(B'\Subset B\), the interior \(H^2\) estimate for the Laplacian, together with \eqref{eq:uniform-local-L2}, gives
\[
\|z_s\|_{H^2(B')}
\le
C_{B',B}
\left(
\|z_s\|_{L^2(B)}
+s^2\|qV_s\|_{L^2(B)}
\right)
\longrightarrow0.
\]
This proves \eqref{eq:zero-frequency-convergence}.

Finally, choose neighborhoods
\[
\partial\Om\Subset\mathcal U'\Subset\mathcal U
\]
such that \(q=1\) almost everywhere on \(\mathcal U\) and
\(\mathcal U\cap\supp f=\varnothing\). On \(\mathcal U\),
\[
(-\Delta+s^2)z_s=-s^2W_{f,c},
\]
and \(W_{f,c}\) is smooth. Iterating the interior estimates for the
constant-coefficient operator \(-\Delta+s^2\), uniformly for
\(0<s\le s_0\), gives
\[
z_s\longrightarrow0
\qquad\text{in }C^\infty(\mathcal U').
\]
Restriction to \(\partial\Om\) proves
\eqref{eq:boundary-zero-frequency}.
\end{proof}

Next we justify the Laplace transform of the distributional boundary trace.
Choose a boundary collar \(\mathcal U'\Subset\mathcal U\) containing
\(\partial\Om\), with \(q=1\) almost everywhere on \(\mathcal U\) and
\(\mathcal U\cap K=\varnothing\). For \(\theta\in C_c^\infty(\R)\), put
\(h_\theta=\int_\R\theta(t)u(t)\,dt\). The constant-coefficient wave equation in this
collar holds for all \(t\in\R\), and hence
\[
\Delta^2 h_\theta=\int_\R\theta^{(4)}(t)u(t)\,dt
\qquad\text{in }\mathcal D'(\mathcal U).
\]
The global bound \eqref{eq:global-wave-bound} and interior elliptic
regularity for \(\Delta^2\) imply
\[
\|h_\theta\|_{H^2(\mathcal U')}
\le C\|f\|_{\mathsf H_{-2}}
\bigl(\|\theta\|_{L^1(\R)}+\|\theta^{(4)}\|_{L^1(\R)}\bigr).
\]
The trace of \(h_\theta\) agrees with the temporal pairing of
\(\Lambda_{f,c}\) with \(\theta\), by the distributional pullback theorem
\cite[Theorem~8.2.4]{HormanderI}. Thus this estimate also bounds
\(\int\theta(t)\Lambda_{f,c}(t,\cdot)\,dt\) in
\(H^{3/2}(\partial\Om)\), uniformly under time translation.
Moreover, finite propagation speed (first for energy data, then by
approximation in \(\mathsf H_{-2}\)) gives
\(\Lambda_{f,c}=0\) for \(0<t<t_*\), for some \(t_*>0\), because
\(K\Subset\Om\). Applying the preceding estimate to smooth cutoffs of
\(e^{-st}\) that vanish before \(t_*/2\) and at large times, and letting
the latter cutoff tend to infinity, shows for each \(s>0\) that the
boundary Laplace transform converges and equals the trace of \(U_s\).
Consequently,
\begin{equation}\label{eq:laplace-trace-identity}
\int_0^\infty e^{-st}\Lambda_{f,c}(t,\cdot)\,dt
=
U_s|_{\partial\Om}
=
sV_s|_{\partial\Om}.
\end{equation}
Here the time integral is understood distributionally; for smooth data, it
agrees with the ordinary integral.

\begin{definition}[Low-frequency boundary coefficient]
\label{def:low-frequency-data}
For a source-speed pair \((f,c)\), define
\begin{equation}\label{eq:low-frequency-boundary-coefficient}
\mathcal L_{f,c}
:=
\lim_{s\downarrow0}
\frac{1}{s}
\int_0^\infty e^{-st}\Lambda_{f,c}(t,\cdot)\,dt
\quad\text{in }C^\infty(\partial\Om).
\end{equation}
Equivalently,
\[
\mathcal L_{f,c}
=
W_{f,c}|_{\partial\Om}.
\]
\end{definition}

\begin{remark}[Observation time]\label{rem:observation-time}
The coefficient \(\mathcal L_{f,c}\) is extracted from the complete
positive-time trace through a Laplace transform. None of the estimates below
is asserted in terms of data on a fixed finite time interval. Such a
formulation would require additional assumptions ensuring that the omitted
late-time contribution is controlled.
\end{remark}

The next lemma converts the Dirichlet trace of a Newtonian potential into
harmonic moment information.

\begin{lemma}[Boundary-to-moment estimate]\label{lem:boundary-to-moment}
Let \(\mu\) be a finite signed or complex Radon measure supported in
\(K\Subset\Om\), and let \(W=\calN\mu\). Then there exists
\(C=C(\Om)>0\) such that
\begin{equation}\label{eq:moment-bound}
\left|\int_\Om\varphi\,d\mu\right|
\le
C\|W\|_{H^{1/2}(\partial\Om)}
 \|\varphi\|_{H^1(\Om)}
\end{equation}
for every harmonic \(\varphi\in H^1(\Om)\).
\end{lemma}

\begin{proof} Since
\[
-\Delta W=\mu
\qquad\text{and}\qquad
-\Delta\varphi=0,
\]
Green's identity formally gives
\begin{equation}\label{eq:green-moment}
\int_\Om\varphi\,d\mu
=
\left\langle
\pnu\varphi,W|_{\partial\Om}
\right\rangle_{\partial\Om}
-
\left\langle
\pnu W,\varphi|_{\partial\Om}
\right\rangle_{\partial\Om},
\end{equation}
where $\langle\cdot,\cdot\rangle_{\partial\Om}$
denotes the duality pairing between \(H^{-1/2}(\partial\Om)\) and
\(H^{1/2}(\partial\Om)\).

This identity is justified by mollifying \(\mu\), applying the weak
Green identity \cite[Lemma~4.3]{McLean2000} to the corresponding smooth
Newtonian potentials, and passing to the limit. Indeed, the mollified
measures may be supported in a fixed compact subset of \(\Om\), and
their potentials converge to \(W\) in \(C^\infty\) near
\(\partial\Om\).

The trace theorem
\cite[Theorem~3.37]{McLean2000} gives
\begin{equation}\label{eq:interior-trace}
\|\varphi|_{\partial\Om}\|_{H^{1/2}(\partial\Om)}
\le
C_\Om\|\varphi\|_{H^1(\Om)}.
\end{equation}
The same theorem provides a bounded right inverse of the trace. Thus,
for every \(\eta\in H^{1/2}(\partial\Om)\), there exists
\(E\eta\in H^1(\Om)\) such that
\[
E\eta|_{\partial\Om}=\eta,
\qquad
\|E\eta\|_{H^1(\Om)}
\le
C_\Om\|\eta\|_{H^{1/2}(\partial\Om)}.
\]
Since \(\varphi\) is harmonic, the definition of its weak normal
derivative \cite[Lemma~4.3]{McLean2000} gives
\[
\left\langle
\pnu\varphi,\eta
\right\rangle_{\partial\Om}
=
\int_\Om\nabla\varphi\cdot\nabla(E\eta)\,dx.
\]
Therefore, by the Cauchy--Schwarz inequality and the boundedness of the
lifting operator,
\[
\begin{aligned}
\left|
\left\langle
\pnu\varphi,\eta
\right\rangle_{\partial\Om}
\right|
&\le
\|\nabla\varphi\|_{L^2(\Om)}
\|\nabla(E\eta)\|_{L^2(\Om)}\\
&\le
C_\Om
\|\varphi\|_{H^1(\Om)}
\|\eta\|_{H^{1/2}(\partial\Om)}.
\end{aligned}
\]
Since this estimate holds for every
\(\eta\in H^{1/2}(\partial\Om)\), the definition of the dual norm yields
\begin{equation}\label{eq:interior-neumann-trace}
\|\pnu\varphi\|_{H^{-1/2}(\partial\Om)}
\le
C_\Om\|\varphi\|_{H^1(\Om)}.
\end{equation}

It remains to estimate the normal derivative of \(W\). Let
\[
\Om^+=\R^3\setminus\overline\Om.
\]
The function \(W\) is the decaying exterior harmonic solution with
Dirichlet data \(W|_{\partial\Om}\). A compactly supported exterior lifting
of these data is provided by the trace right inverse
\cite[Theorem~3.37]{McLean2000} on a bounded boundary collar.
Since \(W\) is harmonic and has finite Dirichlet energy in \(\Om^+\),
the Dirichlet principle bounds its energy by that of this lifting. Thus
\[
\|\nabla W\|_{L^2(\Om^+)}
\le
C_\Om
\|W|_{\partial\Om}\|_{H^{1/2}(\partial\Om)}.
\]
For \(\eta\in H^{1/2}(\partial\Om)\), let \(E^+\eta\) be a compactly
supported exterior lifting satisfying
\[
E^+\eta|_{\partial\Om}=\eta,
\qquad
\|E^+\eta\|_{H^1(\Om^+)}
\le
C_\Om\|\eta\|_{H^{1/2}(\partial\Om)}.
\]
Applying the weak Green identity \cite[Lemma~4.3]{McLean2000}
on a bounded exterior collar containing the support of \(E^+\eta\) gives
\[
\left\langle
\partial_{\nu^+}W,\eta
\right\rangle_{\partial\Om}
=
\int_{\Om^+}
\nabla W\cdot\nabla(E^+\eta)\,dx.
\]
It follows that
\[
\|\partial_{\nu^+}W\|_{H^{-1/2}(\partial\Om)}
\le
C_\Om
\|W|_{\partial\Om}\|_{H^{1/2}(\partial\Om)}.
\]
Since the outward normal \(\nu^+\) of \(\Om^+\) satisfies
\(\nu^+=-\nu\) on \(\partial\Om\), we obtain
\begin{equation}\label{eq:external-dn}
\|\pnu W\|_{H^{-1/2}(\partial\Om)}
\le
C_\Om
\|W|_{\partial\Om}\|_{H^{1/2}(\partial\Om)}.
\end{equation}

Finally, duality in \eqref{eq:green-moment} gives
\[
\begin{aligned}
\left|\int_\Om\varphi\,d\mu\right|
&\le
\|\pnu\varphi\|_{H^{-1/2}(\partial\Om)}
\|W|_{\partial\Om}\|_{H^{1/2}(\partial\Om)}
\\
&\quad+
\|\pnu W\|_{H^{-1/2}(\partial\Om)}
\|\varphi|_{\partial\Om}\|_{H^{1/2}(\partial\Om)}.
\end{aligned}
\]
Applying \eqref{eq:interior-trace},
\eqref{eq:interior-neumann-trace}, and
\eqref{eq:external-dn} proves
\[
\left|\int_\Om\varphi\,d\mu\right|
\le
C_\Om
\|W\|_{H^{1/2}(\partial\Om)}
\|\varphi\|_{H^1(\Om)}.
\]
\end{proof}

\begin{corollary}[Harmonic moments from wave data]
\label{cor:wave-to-moments}
For two source-speed pairs \((f_\ell,c_\ell)\), \(\ell=1,2\), set
\[
q_\ell=c_\ell^{-2},
\qquad
\mu=q_2f_2-q_1f_1,
\qquad
\eps=
\|\mathcal L_{f_2,c_2}-\mathcal L_{f_1,c_1}\|_{H^{1/2}(\partial\Om)}.
\]
Then
\begin{equation}\label{eq:approximate-harmonic-moments}
\left|\int_\Om\varphi\,d\mu\right|
\le
C\eps\|\varphi\|_{H^1(\Om)}
\end{equation}
for every harmonic \(\varphi\in H^1(\Om)\). In particular, if
\[
\Lambda_{f_1,c_1}=\Lambda_{f_2,c_2}
\quad\text{in }\mathcal D'((0,\infty)\times\partial\Om),
\]
then
\[
\int_\Om\varphi\,d\mu=0
\qquad
\text{for every harmonic }\varphi\in H^1(\Om).
\]
\end{corollary}

\begin{proof}
By Proposition~\ref{prop:zero-frequency},
\[
\calN\mu|_{\partial\Om}
=
\mathcal L_{f_2,c_2}-\mathcal L_{f_1,c_1}.
\]
The estimate follows from Lemma~\ref{lem:boundary-to-moment}. Equality of the
time-domain traces gives equality of their Laplace transforms for every
\(s>0\), and hence equality of their low-frequency coefficients.
\end{proof}

\section{Harmonic separation and quantitative interpolation}
\label{sec:harmonic}

The arguments in this section are independent of the wave equation and hold
in every dimension \(n\ge2\). Harmonic functions may be complex valued. We
use \(\zeta\cdot z=\sum_{k=1}^n\zeta_kz_k\) for the complex bilinear
extension of the Euclidean inner product, whereas
\(|\zeta|=(\sum_{k=1}^n|\zeta_k|^2)^{1/2}\) is the Hermitian norm.

If \(\zeta\in\C^n\) satisfies \(\zeta\cdot\zeta=0\), then
\[
\Delta p(\zeta\cdot x)
=
p''(\zeta\cdot x)\,\zeta\cdot\zeta
=
0
\]
for every one-variable polynomial \(p\). These complex harmonic polynomials
will be used both for separation and interpolation.

\begin{lemma}
\label{lem:quantitative-null-vector}
Let \(n\ge2\) and \(M\ge1\). There exists
\(\kappa=\kappa(n,M)>0\) with the following property. For every collection of
nonzero vectors \(v_1,\dots,v_M\in\R^n\), there exists
\(\zeta\in\C^n\) such that
\begin{equation}\label{eq:null-vector-properties}
\zeta\cdot\zeta=0,\qquad |\zeta|=1,
\qquad
|\zeta\cdot v_m|\ge\kappa|v_m|,
\quad 1\le m\le M.
\end{equation}
\end{lemma}

\begin{proof}
Let \(G(2,n)\) denote the Grassmannian of two-dimensional real subspaces of
\(\R^n\), and let \(\pi_P\) denote orthogonal projection onto
\(P\in G(2,n)\). For a tuple
\(\mathbf v=(v_1,\dots,v_M)\in(S^{n-1})^M\), set
\[
F(\mathbf v)
=
\max_{P\in G(2,n)}
\min_{1\le m\le M}|\pi_Pv_m|.
\]
The function \(F\) is continuous. Moreover, \(F(\mathbf v)>0\) for every
tuple. Indeed, for each \(m\), the set
\[
\{P\in G(2,n):P\subset v_m^\perp\}
\]
is closed with empty interior. A finite union of these sets cannot cover the
compact manifold \(G(2,n)\), and one may therefore choose a plane whose
projection of every \(v_m\) is nonzero.

Since \((S^{n-1})^M\) is compact, the positive continuous function \(F\)
has a positive minimum \(a_*=a_*(n,M)\). Apply this conclusion to the
normalized vectors \(v_m/|v_m|\), and choose an orthonormal basis
\(e_1,e_2\) of a plane \(P\) satisfying
\[
|\pi_Pv_m|\ge a_*|v_m|.
\]
Set
\[
\zeta=\frac{e_1+ie_2}{\sqrt2}.
\]
Then \(\zeta\cdot\zeta=0\), \(|\zeta|=1\), and
\[
|\zeta\cdot v_m|^2
=
\frac12\bigl((e_1\cdot v_m)^2+(e_2\cdot v_m)^2\bigr)
=
\frac12|\pi_Pv_m|^2.
\]
Thus \eqref{eq:null-vector-properties} holds with
\(\kappa=a_*/\sqrt2\).
\end{proof}

\begin{lemma}[Finite atomic measures are separated by harmonic functions]
\label{lem:atomic-determination}
Let \(\Om\subset\R^n\), \(n\ge2\), be a domain, and let
\[
\mu=\sum_{\ell=1}^{L}\gamma_\ell\delta_{z_\ell},
\]
where \(z_1,\dots,z_L\in\Om\) are distinct and
\(\gamma_\ell\in\C\). If
\begin{equation}\label{eq:vanishing-harmonic-moments}
\int_\Om\varphi\,d\mu=0
\qquad
\text{for every harmonic }\varphi\text{ in }\Om,
\end{equation}
then \(\mu=0\).
\end{lemma}

\begin{proof}
If \(L=1\), taking \(\varphi\equiv1\) gives \(\gamma_1=0\). Hence assume
\(L\ge2\). Since the points \(z_1,\dots,z_L\) are pairwise distinct, every vector
\(z_\ell-z_m\), with \(\ell\ne m\), is nonzero. Applying
Lemma~\ref{lem:quantitative-null-vector} to this finite collection, choose
a null vector \(\zeta\) such that
\[
\zeta\cdot(z_\ell-z_m)\ne0
\qquad\text{whenever }\ell\ne m.
\]
Therefore, the numbers
\[
\tau_\ell:=\zeta\cdot z_\ell
\]
are pairwise distinct, since
\[
\tau_\ell-\tau_m
=
\zeta\cdot(z_\ell-z_m)\ne0
\qquad\text{for }\ell\ne m.
\] 
For \(k=0,\dots,L-1\), the polynomial
\[
\varphi_k(x)=(\zeta\cdot x)^k
\]
is harmonic. Substitution into \eqref{eq:vanishing-harmonic-moments} gives
\[
\sum_{\ell=1}^{L}\gamma_\ell\tau_\ell^k=0,
\qquad k=0,\dots,L-1.
\]
The coefficient matrix is a Vandermonde matrix with distinct nodes
\(\tau_\ell\), and is therefore invertible. Hence
\(\gamma_1=\cdots=\gamma_L=0\).
\end{proof}

\begin{definition}[Hausdorff distance]\label{def:hausdorff}
For nonempty compact sets \(X,Y\subset\R^n\), their Hausdorff distance is
\[
\dH(X,Y)
=
\max\left\{
\sup_{x\in X}\dist(x,Y),
\sup_{y\in Y}\dist(y,X)
\right\}.
\]
\end{definition}

For the stability analysis, we require a quantitative separation
condition. Let
\[
\mu_1=\sum_{i=1}^{N_1}a_i\delta_{x_i},
\qquad
\mu_2=\sum_{j=1}^{N_2}b_j\delta_{y_j},
\]
where \(x_i,y_j\in K\Subset\Om\). We assume that there exists a prescribed
constant \(\rho>0\) such that distinct points within each configuration
are separated by at least \(\rho\):
\begin{equation}\label{eq:configuration-separation}
|x_i-x_{i'}|\ge\rho
\quad (i\ne i'),
\qquad
|y_j-y_{j'}|\ge\rho
\quad (j\ne j').
\end{equation}
The parameter \(\rho\) is therefore a uniform lower bound on the
within-configuration separation of the source points.

\begin{lemma}[Harmonic interpolation at one atom]
\label{lem:harmonic-interpolation}
Let \(\Om\subset\R^n\) be bounded, let \(K\Subset\Om\), and suppose
\(1\le N_1,N_2\le N_0\) and \eqref{eq:configuration-separation} holds.
Fix \(x_i\) and set
\[
r_i=\dist(x_i,Y),
\qquad
Y=\{y_1,\dots,y_{N_2}\}.
\]
If \(r_i>0\), there exists a harmonic polynomial \(\Phi_i\) such that
\begin{equation}\label{eq:interpolation-values}
\Phi_i(x_i)=1,\qquad
\Phi_i(x_{i'})=0\ (i'\ne i),\qquad
\Phi_i(y_j)=0\ (1\le j\le N_2).
\end{equation}
Moreover,
\begin{equation}\label{eq:coarse-interpolation-bound}
\|\Phi_i\|_{H^1(\Om)}
+
\|\Phi_i\|_{C^1(\overline\Om)}
\le C r_i^{-N_0},
\end{equation}
where \(C\) depends only on
\(\Om,K,n,N_0,\rho\). If, in addition,
\begin{equation}\label{eq:local-distance-regime}
0<r_i\le\frac{\rho}{3},
\end{equation}
then \(\Phi_i\) may be chosen so that the sharper estimate
\begin{equation}\label{eq:sharp-interpolation-bound}
\|\Phi_i\|_{H^1(\Om)}
+
\|\Phi_i\|_{C^1(\overline\Om)}
\le C r_i^{-1}
\end{equation}
holds, with the same parameter dependence.
\end{lemma}

\begin{proof}
Apply Lemma~\ref{lem:quantitative-null-vector} to the family
\[
x_i-x_{i'},\quad i'\ne i,
\qquad
x_i-y_j,\quad 1\le j\le N_2.
\]
There are at most \(2N_0-1\) such vectors. Hence one can choose
\(\zeta\in\C^n\) with \(\zeta\cdot\zeta=0\), \(|\zeta|=1\), and
\begin{equation}\label{eq:projected-difference-bound}
|\zeta\cdot v|\ge\kappa|v|
\end{equation}
for every vector in the family, where
\(\kappa=\kappa(n,N_0)>0\). Define
\begin{equation}\label{eq:interpolation-polynomial}
\Phi_i(x)
=
\prod_{\substack{i'=1\\ i'\ne i}}^{N_1}
\frac{\zeta\cdot(x-x_{i'})}{\zeta\cdot(x_i-x_{i'})}
\prod_{j=1}^{N_2}
\frac{\zeta\cdot(x-y_j)}{\zeta\cdot(x_i-y_j)}.
\end{equation}
This is a polynomial in \(\zeta\cdot x\), and is therefore harmonic.
The interpolation conditions \eqref{eq:interpolation-values} are immediate.

Since \(\Om\) is bounded and all interpolation nodes lie in the fixed compact
set \(K\), the numerator factors and their first derivatives are bounded
uniformly on \(\overline\Om\). The denominators involving
\(x_i-x_{i'}\) are bounded below by \(\kappa\rho\), and those involving
\(x_i-y_j\) are bounded below by \(\kappa r_i\). This proves
\[
\|\Phi_i\|_{H^1(\Om)}
+
\|\Phi_i\|_{C^1(\overline\Om)}
\le Cr_i^{-N_2}.
\]
Because \(N_2\le N_0\), and \(r_i\) is bounded above by a constant depending
on \(\Om\), this implies \eqref{eq:coarse-interpolation-bound}.

Assume now that \(r_i\le\rho/3\), and let \(y_{j_0}\) be a nearest point to
\(x_i\). This point is unique by the separation of \(Y\). For \(j\ne j_0\),
\[
|x_i-y_j|
\ge |y_j-y_{j_0}|-|x_i-y_{j_0}|
\ge \rho-r_i
\ge\frac{2\rho}{3}.
\]
Thus exactly one denominator in the second product in
\eqref{eq:interpolation-polynomial} can be smaller than a fixed multiple of
\(\rho\). Estimate \eqref{eq:projected-difference-bound} then gives
\[
\|\Phi_i\|_{H^1(\Om)}
+
\|\Phi_i\|_{C^1(\overline\Om)}
\le Cr_i^{-1},
\]
which proves \eqref{eq:sharp-interpolation-bound}.
\end{proof}

\begin{remark}
The two estimates in Lemma~\ref{lem:harmonic-interpolation} play different
roles. The coarse estimate is global and places the two configurations within
the separation scale. The sharper estimate is local and is the source of
Lipschitz stability. Using only the first estimate would lead to the weaker
H\"older exponent \(1/N_0\).
\end{remark}

\section{Atomic sources: uniqueness and Lipschitz stability}
\label{sec:atomic}

We now apply the harmonic interpolation results to atomic measures and then
to the wave equation with unknown speed.

\begin{theorem}[Determination of an atomic measure by harmonic moments]
\label{thm:atomic-uniqueness}
Let \(\Om\subset\R^n\), \(n\ge2\), be a domain, and suppose
\[
\mu_1=\sum_{i=1}^{N_1}a_i\delta_{x_i},
\qquad
\mu_2=\sum_{j=1}^{N_2}b_j\delta_{y_j},
\]
where the points in each representation are distinct and all coefficients are
nonzero. If
\begin{equation}\label{eq:equal-harmonic-moments}
\int_\Om\varphi\,d\mu_1
=
\int_\Om\varphi\,d\mu_2
\qquad
\text{for every harmonic }\varphi\text{ in }\Om,
\end{equation}
then \(\mu_1=\mu_2\). In particular,
\[
\{x_1,\dots,x_{N_1}\}
=
\{y_1,\dots,y_{N_2}\},
\]
and, after relabeling, \(a_i=b_i\).
\end{theorem}

\begin{proof}
The finite atomic measure \(\mu=\mu_1-\mu_2\) annihilates every harmonic
function. Combining coefficients at coincident points gives a representation
of \(\mu\) with distinct support points. Lemma~\ref{lem:atomic-determination}
then implies \(\mu=0\).
\end{proof}

\begin{theorem}[Hölder and local Lipschitz stability from approximate harmonic moments]
\label{thm:abstract-lipschitz}
Let \(\Om\subset\R^n\), \(n\ge2\), be a bounded domain, and let
\(K\Subset\Om\). Consider
\[
\mu_1=\sum_{i=1}^{N_1}a_i\delta_{x_i},
\qquad
\mu_2=\sum_{j=1}^{N_2}b_j\delta_{y_j},
\]
where
\[
1\le N_1,N_2\le N_0,
\qquad
x_i,y_j\in K.
\]
Assume the separation condition
\eqref{eq:configuration-separation} and the lower amplitude bounds
\begin{equation}\label{eq:amplitude-lower-bound}
|a_i|\ge m_0,
\qquad
|b_j|\ge m_0.
\end{equation}
Set
\[
X=\{x_1,\dots,x_{N_1}\},
\qquad
Y=\{y_1,\dots,y_{N_2}\},
\qquad
\mu=\mu_2-\mu_1.
\]
Suppose that
\begin{equation}\label{eq:abstract-moment-residual}
\left|\int_\Om\varphi\,d\mu\right|
\le
\eps\|\varphi\|_{H^1(\Om)}
\end{equation}
for every harmonic \(\varphi\in H^1(\Om)\). Then there exist constants
\(C_{\mathrm H},C_{\mathrm L}>0\) and \(\eps_0>0\), depending only on
\[
\Om,\ K,\ n,\ N_0,\ \rho,\ m_0,
\]
such that
\begin{equation}\label{eq:abstract-holder-conclusion}
\dH(X,Y)
\le
C_{\mathrm H}\eps^{1/N_0}
\qquad\text{for every }\eps\ge0,
\end{equation}
and
\begin{equation}\label{eq:abstract-lipschitz-conclusion}
\dH(X,Y)
\le
C_{\mathrm L}\eps
\qquad\text{whenever }0\le\eps\le\eps_0.
\end{equation}
\end{theorem}

\begin{proof}
Fix \(x_i\in X\), and set
\[
r_i=\dist(x_i,Y).
\]
If \(r_i=0\), there is nothing to prove. Otherwise, let \(\Phi_i\) be
the interpolation polynomial from
Lemma~\ref{lem:harmonic-interpolation}. Since \(\Phi_i\) vanishes at
every point of \(Y\) and at every point of \(X\setminus\{x_i\}\), we have
\[
\int_\Om \Phi_i\,d\mu=-a_i.
\]
It follows from \eqref{eq:amplitude-lower-bound},
\eqref{eq:abstract-moment-residual}, and
\eqref{eq:coarse-interpolation-bound} that
\[
m_0
\le |a_i|
\le \eps\|\Phi_i\|_{H^1(\Om)}
\le C\eps r_i^{-N_0}.
\]
Consequently,
\[
r_i\le C\eps^{1/N_0}.
\]
Taking the supremum over \(x_i\in X\), and then repeating the argument
with \(X\) and \(Y\) interchanged, gives
\begin{equation}\label{eq:coarse-distance-estimate}
\dH(X,Y)
\le C_{\mathrm H}\eps^{1/N_0}.
\end{equation}
This proves \eqref{eq:abstract-holder-conclusion} for every
\(\eps\ge0\).

We now obtain the sharper estimate for sufficiently small \(\eps\).
Choose \(\eps_0>0\) such that
\[
C_{\mathrm H}\eps_0^{1/N_0}\le\frac{\rho}{3}.
\]
If \(0\le\eps\le\eps_0\), then
\eqref{eq:coarse-distance-estimate} implies
\[
r_i\le\dH(X,Y)\le\frac{\rho}{3}.
\]
The sharp interpolation estimate
\eqref{eq:sharp-interpolation-bound} therefore applies and gives
\[
m_0
\le |a_i|
\le \eps\|\Phi_i\|_{H^1(\Om)}
\le C\eps r_i^{-1}.
\]
Hence
\[
r_i\le C\eps.
\]
Taking the supremum over \(x_i\in X\), and again interchanging the
roles of \(X\) and \(Y\), yields
\[
\dH(X,Y)\le C_{\mathrm L}\eps,
\]
which proves \eqref{eq:abstract-lipschitz-conclusion}.
\end{proof}

\begin{remark}
The local Lipschitz estimate extends to a global conditional
Lipschitz estimate. Indeed, for $\eps>\eps_0$, since
$X,Y\subset K$, we have
\[
\dH(X,Y)\le \diam(K)
\le \frac{\diam(K)}{\eps_0}\eps.
\]
Combining this with \eqref{eq:abstract-lipschitz-conclusion}
gives
\[
\dH(X,Y)\le C\eps
\qquad\text{for every }\eps\ge0.
\]
The unique matching of individual source points still requires
sufficiently small $\eps$.
\end{remark}

\begin{remark}
A uniform positive lower bound on the nonzero amplitudes is necessary for
support stability: an atom whose coefficient tends to zero becomes
undetectable. Likewise, without a uniform separation condition, two nearby
atoms may merge at the resolution scale of the data, and the interpolation
constants necessarily deteriorate. The constants in
Theorem~\ref{thm:abstract-lipschitz} are therefore conditional on the stated
admissible class.
\end{remark}

\begin{remark}[Cardinality and matching]
For sufficiently small \(\eps\), the preceding proof gives
\[
\dH(X,Y)<\frac{\rho}{3}.
\]
For each \(x_i\in X\), choose a nearest point \(y_{\pi(i)}\in Y\). Then
\[
|x_i-y_{\pi(i)}|
\le \dH(X,Y)
<\frac{\rho}{3}.
\]
The map \(i\mapsto\pi(i)\) is injective. Indeed, if two distinct points
\(x_i\) and \(x_{i'}\) had the same nearest point \(y_{\pi(i)}\), then
\[
|x_i-x_{i'}|
\le
|x_i-y_{\pi(i)}|
+
|y_{\pi(i)}-x_{i'}|
<
\frac{2\rho}{3},
\]
contradicting the separation condition
\(|x_i-x_{i'}|\ge\rho\). Hence \(N_1\le N_2\). Interchanging the roles
of \(X\) and \(Y\) gives \(N_2\le N_1\), and therefore
\[
N_1=N_2.
\]
Thus the nearest-point map is a bijection. After relabeling the points of
\(Y\) accordingly, the local Lipschitz estimate yields
\[
|x_i-y_i|
\le
\dH(X,Y)
\le
C_{\mathrm L}\eps,
\qquad i=1,\dots,N_1.
\]
\end{remark}

We return to the wave equation. Let
\begin{equation}\label{eq:two-atomic-sources}
f_1=\sum_{i=1}^{N_1}\alpha_i\delta_{x_i},
\qquad
f_2=\sum_{j=1}^{N_2}\beta_j\delta_{y_j},
\end{equation}
and let \(q_\ell=c_\ell^{-2}\). The weighted measures are
\[
q_1f_1
=
\sum_{i=1}^{N_1}\alpha_iq_1(x_i)\delta_{x_i},
\qquad
q_2f_2
=
\sum_{j=1}^{N_2}\beta_jq_2(y_j)\delta_{y_j}.
\]

\begin{theorem}[Unknown-speed localization and weighted-amplitude stability]
\label{thm:wave-atomic-localization}
Assume the hypotheses of Section~\ref{sec:low-frequency} for
\((f_\ell,c_\ell)\), \(\ell=1,2\), and suppose that the sources are given
by \eqref{eq:two-atomic-sources}, with nonzero coefficients.

\begin{enumerate}[label=\textup{(\roman*)}]
\item If
\[
\Lambda_{f_1,c_1}
=
\Lambda_{f_2,c_2}
\quad\text{in }
\mathcal D'((0,\infty)\times\partial\Om),
\]
then
\[
\supp f_1=\supp f_2,
\qquad
q_1f_1=q_2f_2.
\]
Thus the source locations and the weighted amplitudes are uniquely
determined.

\item Fix \(N_0\), \(\rho\), and \(m_0>0\). Suppose that
\(N_1,N_2\le N_0\), the point sets satisfy
\eqref{eq:configuration-separation}, and
\begin{equation}\label{eq:weighted-amplitude-lower}
|\alpha_iq_1(x_i)|\ge m_0,
\qquad
|\beta_jq_2(y_j)|\ge m_0.
\end{equation}
Set
\[
\delta
=
\|\mathcal L_{f_2,c_2}
-\mathcal L_{f_1,c_1}\|_{H^{1/2}(\partial\Om)}.
\]
There exist \(C_{\mathrm H},C_{\mathrm L}>0\) and \(\eps_0>0\),
depending only on
\[
\Om,\ K,\ N_0,\ \rho,\ m_0,
\]
such that
\begin{equation}\label{eq:wave-atomic-holder-stability}
\dH(\supp f_1,\supp f_2)
\le
C_{\mathrm H}\delta^{1/N_0}
\end{equation}
for every \(\delta\ge0\), and
\begin{equation}\label{eq:wave-atomic-lipschitz-stability}
\dH(\supp f_1,\supp f_2)
\le
C_{\mathrm L}\delta
\end{equation}
whenever \(0\le\delta\le\eps_0\).

For \(\eps_0\) sufficiently small, the latter condition also implies
\(N_1=N_2=:N\), and the points may be relabeled so that
\begin{equation}\label{eq:wave-atomic-pointwise-matching}
|x_i-y_i|
\le
C_{\mathrm L}\delta,
\qquad i=1,\dots,N.
\end{equation}
If, in addition, the weighted amplitudes satisfy
\begin{equation}\label{eq:weighted-amplitude-upper}
|\alpha_iq_1(x_i)|\le M_0,
\qquad
|\beta_jq_2(y_j)|\le M_0,
\end{equation}
then
\begin{equation}\label{eq:wave-atomic-weighted-amplitude-stability}
\max_{1\le i\le N}
\left(
|x_i-y_i|
+
|\alpha_iq_1(x_i)-\beta_iq_2(y_i)|
\right)
\le
C\delta,
\end{equation}
where \(C\) may also depend on \(M_0\).
\end{enumerate}
\end{theorem}

\begin{proof}
Set
\[
\mu_1=q_1f_1
=
\sum_{i=1}^{N_1}\alpha_iq_1(x_i)\delta_{x_i},
\qquad
\mu_2=q_2f_2
=
\sum_{j=1}^{N_2}\beta_jq_2(y_j)\delta_{y_j},
\]
and let
\[
\mu=\mu_2-\mu_1.
\]

If the time-domain boundary measurements agree,
Corollary~\ref{cor:wave-to-moments} shows that
\[
\int_\Om\varphi\,d\mu=0
\]
for every harmonic \(\varphi\in H^1(\Om)\). Since the proof of
Theorem~\ref{thm:atomic-uniqueness} uses only harmonic polynomials,
which belong to \(H^1(\Om)\) because \(\Om\) is bounded, it follows
that \(\mu=0\), or equivalently,
\[
q_1f_1=q_2f_2.
\]
Since \(q_\ell>0\) and all source coefficients are nonzero, multiplication
by \(q_\ell\) neither creates nor removes support points. Therefore
\[
\supp f_1=\supp f_2,
\]
which proves part~\textup{(i)}.

For part~\textup{(ii)}, Corollary~\ref{cor:wave-to-moments} gives
\[
\left|\int_\Om\varphi\,d\mu\right|
\le
C_\Om\delta\,
\|\varphi\|_{H^1(\Om)}
\]
for every harmonic \(\varphi\in H^1(\Om)\). Thus
\eqref{eq:abstract-moment-residual} holds with
\(\eps=C_\Om\delta\). The weighted amplitudes of \(\mu_1\) and
\(\mu_2\) satisfy the lower bound
\eqref{eq:amplitude-lower-bound} by
\eqref{eq:weighted-amplitude-lower}. Applying
Theorem~\ref{thm:abstract-lipschitz} and absorbing \(C_\Om\) into the
constants gives \eqref{eq:wave-atomic-holder-stability} and
\eqref{eq:wave-atomic-lipschitz-stability}.

After decreasing \(\eps_0\), if necessary, we may assume that
\[
\dH(\supp f_1,\supp f_2)<\frac{\rho}{3}
\qquad\text{whenever}\qquad
\delta\le\eps_0.
\]
The separation condition then implies, by the nearest-point matching
argument, that \(N_1=N_2=:N\). After relabeling,
\eqref{eq:wave-atomic-lipschitz-stability} gives
\[
|x_i-y_i|
\le C_{\mathrm L}\delta,
\qquad i=1,\dots,N.
\]

It remains to estimate the weighted amplitudes. Write
\[
a_i=\alpha_iq_1(x_i),
\qquad
b_i=\beta_iq_2(y_i).
\]
Fix a ball \(B\) containing \(\overline\Om\). For each \(i\), the quantitative
harmonic interpolation construction underlying
Lemma~\ref{lem:harmonic-interpolation} provides a harmonic polynomial
\(\Psi_i\) such that
\[
\Psi_i(x_j)=\delta_{ij},
\qquad
\|\Psi_i\|_{H^1(\Om)}
+
\|\nabla\Psi_i\|_{L^\infty(B)}
\le C,
\]
where \(C\) is uniform over all configurations satisfying the prescribed
separation condition. We have
\[
\int_\Om\Psi_i\,d\mu
=
\sum_{j=1}^{N}b_j\Psi_i(y_j)-a_i.
\]
Using \(\Psi_i(x_j)=\delta_{ij}\), it follows that
\[
b_i-a_i
=
\int_\Om\Psi_i\,d\mu
-
\sum_{j=1}^{N}
b_j\bigl(\Psi_i(y_j)-\Psi_i(x_j)\bigr).
\]
Therefore, by the moment estimate, the upper amplitude bound
\eqref{eq:weighted-amplitude-upper}, and the mean-value theorem,
\[
\begin{aligned}
|b_i-a_i|
&\le
C_\Om\delta\|\Psi_i\|_{H^1(\Om)}
+
\sum_{j=1}^{N}
|b_j|\,
|\Psi_i(y_j)-\Psi_i(x_j)| \\
&\le
C\delta
+
CM_0\sum_{j=1}^{N}|x_j-y_j| \\
&\le
C\delta.
\end{aligned}
\]
Combining this estimate with
\eqref{eq:wave-atomic-pointwise-matching} proves
\eqref{eq:wave-atomic-weighted-amplitude-stability}.
\end{proof}

\begin{remark}
For atomic sources, the low-frequency coefficient determines the weighted
source \(qf\). After matching the common support points, this means that
\[
\alpha_iq_1(x_i)
=
\beta_iq_2(x_i),
\qquad\text{or equivalently}\qquad
\alpha_i c_1(x_i)^{-2}
=
\beta_i c_2(x_i)^{-2}.
\]
The coefficient alone cannot separate the physical amplitudes from the
unknown speed factors at these points.
\end{remark}

\section{Mollified point sources}
\label{sec:mollified}

We next replace each point mass by a small smooth profile. In contrast to the
atomic case, exact recovery from harmonic moments cannot hold for arbitrary
smooth sources. The points \(x_i\) and \(y_j\) below are translation
parameters of the profiles; for a nonradial profile they need not coincide
with a centroid or a point of maximum amplitude. The appropriate conclusion
is stable localization of these translation parameters, up to an error
proportional to the size of the profiles.

Let \(\Om\subset\R^n\), \(n\ge2\), be bounded, and let \(K\Subset\Om\).
For \(\eta>0\) and \(\chi\in C_c^\infty(B(0,1))\), set
\[
\chi_\eta(x)=\eta^{-n}\chi(x/\eta).
\]
Thus, for example, \(\chi_{i,\eta_i}(x-x_i)\) denotes
\(\eta_i^{-n}\chi_i((x-x_i)/\eta_i)\).

\begin{assumption}[Admissible separated blobs]
\label{ass:admissible-blobs}
The two sources have the form
\begin{equation}\label{eq:mollified-sources}
f_1(x)
=
\sum_{i=1}^{N_1}\alpha_i
\chi_{i,\eta_i}(x-x_i),
\qquad
f_2(x)
=
\sum_{j=1}^{N_2}\beta_j
\psi_{j,\sigma_j}(x-y_j),
\end{equation}
where
\[
1\le N_1,N_2\le N_0,
\qquad x_i,y_j\in K,
\]
and
\[
\chi_i,\psi_j\in C_c^\infty(B(0,1)),
\qquad
\chi_i,\psi_j\ge0,
\qquad
\int_{\R^n}\chi_i\,dx
=
\int_{\R^n}\psi_j\,dx
=1.
\]
Set
\[
\eta_*=
\max\{\eta_1,\dots,\eta_{N_1},
       \sigma_1,\dots,\sigma_{N_2}\}.
\]
We assume
\begin{equation}\label{eq:blob-radius-interior}
\eta_*
\le
\frac12\dist(K,\partial\Om),
\end{equation}
the separation of the translation parameters
\begin{equation}\label{eq:blob-separation}
|x_i-x_{i'}|\ge\rho\quad(i\ne i'),
\qquad
|y_j-y_{j'}|\ge\rho\quad(j\ne j'),
\qquad
\eta_*\le\frac{\rho}{8},
\end{equation}
the amplitude bounds
\begin{equation}\label{eq:blob-amplitude-bounds}
0<\alpha_0\le|\alpha_i|,|\beta_j|\le A_0,
\end{equation}
and the speed bounds
\begin{equation}\label{eq:blob-speed-bounds}
0<c_-\le c_\ell(x)\le c_+<\infty,
\qquad \ell=1,2.
\end{equation}
\end{assumption}

Let
\[
X=\{x_1,\dots,x_{N_1}\},
\qquad
Y=\{y_1,\dots,y_{N_2}\},
\]
and define
\begin{equation}\label{eq:blob-weighted-difference}
d\mu
=
\left(c_2^{-2}f_2-c_1^{-2}f_1\right)\,dx.
\end{equation}

\begin{theorem}[Lipschitz localization of blob translation parameters]
\label{thm:blob-stability}
Assume \ref{ass:admissible-blobs}. Suppose that, for some \(\eps\ge0\),
\begin{equation}\label{eq:blob-moment-residual}
\left|\int_\Om\varphi\,d\mu\right|
\le
\eps\|\varphi\|_{H^1(\Om)}
\end{equation}
for every harmonic \(\varphi\in H^1(\Om)\). Then there exist
\(C>0\) and \(\delta_0>0\), depending only on
\[
\Om,\ K,\ n,\ N_0,\ \rho,\ \alpha_0,\ A_0,\ c_-,\ c_+,
\]
such that, if
\[
\eps+\eta_*\le\delta_0,
\]
then
\begin{equation}\label{eq:blob-center-lipschitz}
\dH(X,Y)\le C(\eps+\eta_*).
\end{equation}
Moreover,
\begin{equation}\label{eq:blob-support-lipschitz}
\dH(\supp f_1,\supp f_2)
\le C(\eps+\eta_*).
\end{equation}
\end{theorem}

\begin{proof}
Write
\[
g_i^{(1)}(x)
=
\alpha_i c_1(x)^{-2}\chi_{i,\eta_i}(x-x_i),
\qquad
g_j^{(2)}(x)
=
\beta_j c_2(x)^{-2}\psi_{j,\sigma_j}(x-y_j),
\]
and define their weighted masses
\[
m_i^{(1)}=\int_\Om g_i^{(1)}(x)\,dx,
\qquad
m_j^{(2)}=\int_\Om g_j^{(2)}(x)\,dx.
\]
Since the profiles are nonnegative and normalized,
\[
m_i^{(1)}
=
\alpha_i\int_\Om c_1(x)^{-2}
\chi_{i,\eta_i}(x-x_i)\,dx,
\]
and the analogous identity holds for \(m_j^{(2)}\). Therefore
\begin{equation}\label{eq:blob-mass-lower}
|m_i^{(1)}|,
|m_j^{(2)}|
\ge
\alpha_0c_+^{-2}
=:m_*.
\end{equation}
Similarly,
\begin{equation}\label{eq:blob-l1-upper}
\|g_i^{(1)}\|_{L^1(\Om)},
\|g_j^{(2)}\|_{L^1(\Om)}
\le
A_0c_-^{-2}
=:M_*.
\end{equation}

Fix \(x_i\in X\) and set \(r_i=\dist(x_i,Y)\). If \(r_i=0\), there is
nothing to prove. Let \(\Phi_i\) be the harmonic interpolation polynomial
from Lemma~\ref{lem:harmonic-interpolation}. Since
\(\Phi_i(x_i)=1\) and \(\Phi_i\) vanishes at all the remaining centers,
\[
\begin{aligned}
m_i^{(1)}
={}&
-\int_\Om\Phi_i\,d\mu\\
&+
\sum_{j=1}^{N_2}
\int_\Om
\bigl(\Phi_i(x)-\Phi_i(y_j)\bigr)g_j^{(2)}(x)\,dx\\
&-
\sum_{\substack{\ell=1\\ \ell\ne i}}^{N_1}
\int_\Om
\bigl(\Phi_i(x)-\Phi_i(x_\ell)\bigr)g_\ell^{(1)}(x)\,dx\\
&-
\int_\Om
\bigl(\Phi_i(x)-\Phi_i(x_i)\bigr)g_i^{(1)}(x)\,dx.
\end{aligned}
\]
Every weighted blob is supported in the ball of its translation parameter with
radius at
most \(\eta_*\). By the mean value theorem and
\eqref{eq:blob-l1-upper},
\[
\left|
\int_\Om
\bigl(\Phi_i(x)-\Phi_i(z)\bigr)g(x)\,dx
\right|
\le
M_*\eta_*\|\Phi_i\|_{C^1(\overline\Om)}
\]
for each weighted blob \(g=g_l^{(1)}\) or \(g=g_j^{(2)}\),
with corresponding reference point \(z=x_l\) or \(z=y_j\). Using
\eqref{eq:blob-moment-residual} and the fact that there are at most
\(2N_0\) blobs, we obtain
\begin{equation}\label{eq:blob-main-estimate}
m_*
\le
C(\eps+\eta_*)
\left(
\|\Phi_i\|_{H^1(\Om)}
+
\|\Phi_i\|_{C^1(\overline\Om)}
\right).
\end{equation}

The coarse interpolation bound \eqref{eq:coarse-interpolation-bound} gives
\[
r_i
\le
C(\eps+\eta_*)^{1/N_0}.
\]
Repeating the argument after interchanging \(X\) and \(Y\) yields
\begin{equation}\label{eq:blob-coarse-hausdorff}
\dH(X,Y)
\le
C(\eps+\eta_*)^{1/N_0}.
\end{equation}
Choose \(\delta_0\) sufficiently small that the right-hand side is less than
\(\rho/3\). Then \(r_i\le\rho/3\), and the sharp interpolation estimate
\eqref{eq:sharp-interpolation-bound}, together with
\eqref{eq:blob-main-estimate}, gives
\[
m_*\le C(\eps+\eta_*)r_i^{-1}.
\]
Thus
\[
r_i\le C(\eps+\eta_*).
\]
Taking the supremum over \(x_i\), and repeating the argument with the two
configurations interchanged, proves \eqref{eq:blob-center-lipschitz}.

It remains to estimate the Hausdorff distance between the supports of the two sources using the bound on the Hausdorff distance between their translation-parameter sets. Because of
\(\eta_*\le\rho/8\), the blobs belonging to each source have disjoint
supports, and no cancellation between different terms can remove a blob from
the support. Every point of \(\supp f_1\) lies within \(\eta_*\) of \(X\).
Conversely, each normalized profile has nonempty support in \(B(0,1)\), so
every translation parameter in \(X\) lies within \(\eta_*\) of \(\supp f_1\). Hence
\[
\dH(\supp f_1,X)\le\eta_*.
\]
Similarly,
\[
\dH(\supp f_2,Y)\le\eta_*.
\]
The triangle inequality for Hausdorff distance and
\eqref{eq:blob-center-lipschitz} yield
\[
\dH(\supp f_1,\supp f_2)
\le
2\eta_*+\dH(X,Y)
\le C(\eps+\eta_*),
\]
which proves \eqref{eq:blob-support-lipschitz}.
\end{proof}

\begin{corollary}[Blob translation-parameter localization from unknown-speed wave data]
\label{cor:blob-wave-stability}
Assume \(n=3\), the wave-speed hypotheses of Section~\ref{sec:low-frequency},
and Assumption~\ref{ass:admissible-blobs}. Set
\[
\eps_{\partial}
=
\|\mathcal L_{f_2,c_2}-\mathcal L_{f_1,c_1}\|_{H^{1/2}(\partial\Om)}.
\]
There exist \(C,\delta_0>0\), depending only on
\(\Om,K,n,N_0,\rho,\alpha_0,A_0,c_-,c_+\), such that
\[
\dH(X,Y)
+
\dH(\supp f_1,\supp f_2)
\le
C(\eps_{\partial}+\eta_*)
\]
whenever \(\eps_{\partial}+\eta_*\le\delta_0\).
\end{corollary}

\begin{proof}
By \eqref{eq:blob-radius-interior}, both sources are supported in the fixed
compact set
\[
K'=\{x\in\R^3:\dist(x,K)\le\tfrac12\dist(K,\partial\Om)\}
\Subset\Om.
\]
Corollary~\ref{cor:wave-to-moments}, applied with \(K'\), gives
\eqref{eq:blob-moment-residual} with
\(\eps\le C\eps_{\partial}\). Apply
Theorem~\ref{thm:blob-stability}.
\end{proof}

\begin{remark}[Exact data and finite source size]
Even if the boundary data agree exactly, Theorem~\ref{thm:blob-stability}
gives an error of order \(\eta_*\). This is the natural resolution scale of
the present admissible source class; no sharper conclusion is asserted here
for general nonradial profiles.
\end{remark}

\begin{remark}[Weighted amplitudes]
In the small-error regime of Theorem~\ref{thm:blob-stability}, the
nearest-point argument also gives \(N_1=N_2=:N\). After relabeling the
blobs according to this matching, the interpolation argument can be
extended to show stability of the weighted blob masses
\[
m_i^{(1)}
=
\alpha_i\int_\Om q_1(x)
\chi_{i,\eta_i}(x-x_i)\,dx,
\qquad
m_i^{(2)}
=
\beta_i\int_\Om q_2(x)
\psi_{i,\sigma_i}(x-y_i)\,dx.
\]
More precisely,
\[
|m_i^{(1)}-m_i^{(2)}|
\le
C(\eps+\eta_*),
\qquad i=1,\dots,N.
\]
For merely bounded measurable \(q_\ell\), these averaged quantities are
the natural weighted amplitudes of the blobs. If \(q_1\) and \(q_2\) are
Lipschitz on a fixed neighborhood containing all blob supports, then
\[
|m_i^{(1)}-\alpha_iq_1(x_i)|
+
|m_i^{(2)}-\beta_iq_2(y_i)|
\le
C\eta_*,
\]
and consequently
\[
|\alpha_iq_1(x_i)-\beta_iq_2(y_i)|
\le
C(\eps+\eta_*).
\]
Thus pointwise weighted-amplitude stability is also available under this
additional local regularity assumption on the speeds.
\end{remark}

\section*{Acknowledgments} The author is partially supported by NSF grant DMS-1953620.
	
\subsection*{Use of artificial intelligence}
	
	OpenAI's ChatGPT 5.6 Sol was used during the preparation of this work to assist with editorial refinement and identification of relevant literature. The main theoretical results were written by the author, who also verified and edited all output. The author assumes responsibility for all content.

\end{document}